\documentclass{article}
\usepackage{amsfonts}
\usepackage{amssymb}
\usepackage{amsmath}
\usepackage{hyperref}
\usepackage{geometry}

\newtheorem{theorem}{Theorem}

\newtheorem{corollary}[theorem]{Corollary}

\newtheorem{lemma}[theorem]{Lemma}

\newtheorem{problem}[theorem]{Problem}
\newtheorem{proposition}[theorem]{Proposition}

\newenvironment{proof}[1][Proof]{\noindent\textbf{#1.} }{\ \rule{0.5em}{0.5em}}
\input{tcilatex}
\begin{document}

\title{{\Huge Compactness of the space of left orders}\bigskip \\
{\Large 
\centerline{Dedicated to L. H. Kauffman for his 60th birthday}%
}}
\author{M. A. Dabkowska \\
Department of Mathematical Sciences\\
University of Texas at Dallas\\
Richardson, TX 75083\\
madab@utdallas.edu \and M. K. Dabkowski \\
Department of Mathematical Sciences\\
University of Texas at Dallas\\
Richardson, TX 75083\\
mdab@utdallas.edu \and V. S. Harizanov \\
Department of Mathematics\\
George Washington University\\
Washington, DC 20052\\
harizanv@gwu.edu \and J. H. Przytycki \\
Department of Mathematics\\
George Washington University\\
Washington, DC 20052\\
przytyck@gwu.edu \and M. A. Veve \\
Department of Mathematics\\
George Washington University\\
Washington, DC 20052\\
veve@gwu.edu }
\date{}
\maketitle

\begin{abstract}
A left order on a magma (e.g., semigroup) is a total order of its elements
that is left invariant under the magma operation. A natural topology can be
introduced on the set of all left orders of an arbitrary magma. We prove
that this topological space is compact. Interesting examples of
nonassociative magmas, whose spaces of right orders we analyze, come from
knot theory and are called quandles. Our main result establishes an
interesting connection between topological properties of the space of left
orders on a group, and the classical algebraic result by Conrad \cite{Conrad}
and \L o\'{s} \cite{Los} concerning the existence of left orders.\bigskip

\noindent 
\it{Keywords}: \textup{magma, order, quandle, topology on orders, Cantor cube}%
\bigskip

\noindent 
\textup{Mathematics Subject Classification 2000: Primary 57M99; Secondary 17A99, 03G05, 54D30}%
\end{abstract}

\section{Introduction}

\noindent In recent years, the theory of orders on groups has become an
important tool in understanding the geometric properties of $3$-dimensional
manifolds (see \cite{Rolfsen}). For a semigroup $G$, A. S. Sikora \cite%
{Sikora} defined a natural topology on the space of its left orders $LO(G)$.
He showed that if $G$ is countable, then this space is compact, metrizable,
and totally disconnected. For an arbitrary (not necessarily associative)
magma $\mathcal{M}$, we analyze the space of its left (right) orders $LO(%
\mathcal{M})$ ($RO(\mathcal{M)}$). We show that $LO(\mathcal{M})$ is a
compact topological space that can be embedded into the Cantor cube $%
\{0,1\}^{\mathfrak{m}}$, where $\mathfrak{m=}\left\vert \mathcal{M}%
\right\vert $ is the cardinality of $\mathcal{M}$. In the case of a group $G,
$ this result is related to the classical theorem of P. Conrad in \cite%
{Conrad}. Conrad's theorem provides a necessary and sufficient algebraic
condition (in terms of semigroups) for a partial left order to extend to a
linear (total) left order on $G$. This connection relies on Alexander's
Subbase Theorem \cite{Alex}. Interesting examples of nonassociative magma,
whose spaces of right orders we analyze, come from knot theory and are
called quandles \cite{Joyce}.

Let $\mathcal{M}$ be a magma\footnote{%
The term magma was used by J-P.~Serre \cite{Ser} and Bourbaki \cite{Bou}.},
that is, $\mathcal{M}$ is a set with a binary operation $\cdot :\mathcal{%
M\times M\rightarrow M}$. We consider all strict linear orders $\mathcal{R}$
on $\mathcal{M}$ which are left invariant under the magma operation:
\begin{equation*}
(\forall a,b,c\in \mathcal{M)[(}a,b\mathcal{)\in R}\Rightarrow (ca,cb)\in 
\mathcal{R}].
\end{equation*}%
As usual, $LO(\mathcal{M})$ denotes the set of all left invariant strict
orders on $\mathcal{M}$, $RO(\mathcal{M})$ the set of all right invariant
strict orders, and $BiO(\mathcal{M})$ the set of all bi-orders. Clearly, $%
BiO(\mathcal{M})=LO(\mathcal{M})\cap RO(\mathcal{M})$.

We define a topology on $LO(\mathcal{M})$ by choosing as a subbasis the
collection 
$\mathcal{S=}\{S_{(a,b)}\}_{(a,b)\in (\mathcal{M\times M)}\backslash \Delta
}$, where $S_{(a,b)}=\{\mathcal{R}\in LO(\mathcal{M})|$ $(a,b)\in \mathcal{R}%
\}$ and $\Delta =\{(a,a)$ $|$ $a\in \mathcal{M}\}$. Recall that a
topological space is zero-dimensional if it is a $T_{1}$-space with a clopen
(closed and open) basis. The following proposition follows directly from the
definition.

\begin{proposition}
\label{Zero_dim}The space $LO(\mathcal{M})$ is zero-dimensional.
\end{proposition}

\section{Main result}

We show that for any magma $\mathcal{M}$, the space of left invariant orders 
$LO(\mathcal{M})$ is a compact topological space. Recall that the \emph{%
weight} of a topological space $X$ is the minimal cardinality $\mathfrak{%
\kappa }$ of a basis for the topology on $X$. For technical convenience the
weight is defined to be $\aleph _{0}$ when the minimal basis is finite.
Every zero-dimensional space of weight $\mathfrak{m}$ has a clopen subbasis
of cardinality $\mathfrak{m}$. We use the result by Vedenissoff \cite{Ved}
that if $X$ is a zero-dimensional space of weight $\mathfrak{m}$, then $X$
can be embedded into the Cantor cube $\{0,1\}^{\mathfrak{m}}$. Such a
homeomorphic embedding is defined as follows. Let $\{U_{\alpha }\}_{\alpha
\in \Gamma }$ be a clopen subbasis of the cardinality $\mathfrak{m}$. For
every $\alpha \in \Gamma $, define a mapping $\psi _{\alpha }:X\rightarrow
\{0,1\}_{\alpha },$ where $\{0,1\}_{\alpha }=\{0,1\}$ for $\alpha \in
\Gamma $, as follows:%
\begin{equation*}
\psi _{\alpha }(x)=\left\{ 
\begin{array}{ccc}
1 & if & x\in U_{\alpha } \\ 
0 & if & x\notin U_{\alpha }%
\end{array}%
\right. 
\end{equation*}%
The mapping $\psi =\prod_{\alpha \in \Gamma }\psi _{\alpha }$ is a
homeomorphic embedding of $X$ into $\{0,1\}^{\mathfrak{m}}=$ $\prod_{\alpha
\in \Gamma }\{0,1\}_{\alpha }$ (see \cite{Eng}; Theorem 6.2.16).

\begin{theorem}
\label{Embedding}Let $\mathcal{M}$ be a magma with $\left\vert \mathcal{M}%
\right\vert =\mathfrak{m\geq \aleph }_{0}$. Then $LO(\mathcal{M})$ is a
closed subspace of the Cantor cube $\{0,1\}^{\mathfrak{m}}$. In particular, $%
LO(\mathcal{M})$ is a compact space.
\end{theorem}

\begin{proof}%
By Proposition \ref{Zero_dim}, $LO(\mathcal{M})$ is a zero-dimensional space
with a clopen subbasis

\noindent $\mathcal{S=}\{S_{(a,b)}\}_{(a,b)\in (\mathcal{M\times M)}%
\backslash \Delta }$ of weight $\left\vert \mathcal{S}\right\vert \leq 
\mathfrak{m}$. By Vedenissoff's theorem \cite{Ved}, $LO(\mathcal{M})$ can be
homeomorphically embedded into the Cantor cube $\{0,1\}^{\mathfrak{m}}$.
Since the Cantor cube $\{0,1\}^{\mathfrak{m}}$ is compact, it suffices to
show that the image of the embedding $\psi :LO(\mathcal{M})\rightarrow
\{0,1\}^{\mathfrak{m}}$ is a closed subspace of $\{0,1\}^{\mathfrak{m}}$
(or, equivalently, $\{0,1\}^{\mathfrak{m}}\backslash \psi (LO(\mathcal{M}))$
is open). Each $\mathcal{R\in }LO\mathcal{(\mathcal{M})}$ can be viewed as
the characteristic function $\chi _{\mathcal{R}}$ of $\mathcal{R}\subset (%
\mathcal{M\times M)}\backslash \Delta $:%
\begin{equation*}
\chi _{\mathcal{R}}(a,b)=\left\{ 
\begin{array}{ccc}
1 & if & (a,b)\in \mathcal{R} \\ 
0 & if & (a,b)\notin \mathcal{R}%
\end{array}%
\right. 
\end{equation*}%
In particular,%
\begin{equation*}
\mathbf{(i)}\text{ }(\forall c\in \mathcal{M})[\chi _{\mathcal{R}%
}(a,b)=1\Rightarrow \chi _{\mathcal{R}}(ca,cb)=1]
\end{equation*}%
and%
\begin{equation*}
\mathbf{(ii)}\text{ }(\forall c\in \mathcal{M})[(\chi _{\mathcal{R}%
}(a,b)=1)\wedge (\chi _{\mathcal{R}}(b,c)=1)\Rightarrow (\chi _{\mathcal{R}%
}(a,c)=1)],
\end{equation*}%
and%
\begin{equation*}
\mathbf{(iii)}\text{ }(\forall a,b\in \mathcal{M})[\chi _{\mathcal{R}%
}(a,b)\neq \chi _{\mathcal{R}}(b,a)]
\end{equation*}%
Let $f\in \{0,1\}^{\mathfrak{m}}\backslash \psi (LO(\mathcal{M}))$. This
means that $f:(\mathcal{M\times M)}\backslash \Delta \rightarrow \{0,1\}$
is not a characteristic function for any strict left order on $\mathcal{M}$.
Therefore, $f$ fails to satisfy at least one of the conditions $\mathbf{(i)}$%
, $\mathbf{(ii)}$ and $\mathbf{(iii)}$ given above.

First we assume that $f$ fails to satisfy condition $\mathbf{(i)}$. Then
there are $a,b,c\in \mathcal{M}$ such that%
\begin{equation*}
f(a,b)=1\text{ and }f(ca,cb)=0\text{, or}
\end{equation*}%
\begin{equation*}
f(a,b)=1\text{ and }ca=cb.
\end{equation*}

\begin{itemize}
\item If $f(a,b)=1$ and $ca=cb$, where $%
a\neq b$, then $LO(\mathcal{M})=\emptyset $.

\item If $f(a,b)=1$ and $f(ca,cb)=0$, then the set $U_{f}$ defined by 
\begin{equation*}
U_{f}=\{g\in \{0,1\}^{\mathfrak{m}}|\text{ }g(a,b)=1\text{ and }g(ca,cb)=0\}
\end{equation*}%
is open in the Cantor cube. Moreover, $f\in U_{f}$ and $U_{f}\subset
\{0,1\}^{\mathfrak{m}}\backslash \psi (LO(\mathcal{M}))$.
\end{itemize}

\noindent Now suppose that $f$ fails to satisfy $\mathbf{(ii)}$ or $\mathbf{%
(iii)}$. This means that:

\begin{description}
\item[(a)] There are $a,b,c\in \mathcal{M}$ such that $f(a,b)=f(b,c)=1$ but $%
f(a,c)=0,$ or

\item[(b)] There are $a,b\in \mathcal{M}$ such that $f(a,b)=f(b,a)$.
\end{description}

\noindent In case $(a)$ we can argue as above. Namely, the set 
\begin{equation*}
U_{f}=\{g\in \{0,1\}^{\mathfrak{m}}|\text{ }g(a,b)=g(b,c)=1\text{ and }%
g(a,c)=0\}
\end{equation*}%
is an open neighborhood disjoint from the image of $\psi $. A similar
argument applies to the case $(b)$. Hence $LO(\mathcal{M})$ is a closed
subspace of the Cantor cube $\{0,1\}^{\mathfrak{m}}$. Therefore, $LO(%
\mathcal{M})$ is compact.%
\end{proof}%

\begin{corollary}
The space of bi-orders $BiO(\mathcal{M})$ is compact.
\end{corollary}

\begin{proof}%
Similarly to the proof of Theorem \ref{Embedding} we show that $RO(\mathcal{M%
})$ is compact. Therefore, $BiO(\mathcal{M})=RO(\mathcal{M})\cap LO(\mathcal{%
M})$ is also compact.%
\end{proof}%

\begin{corollary}
\label{Metr}If $\mathcal{M}$ is a countable magma, then $LO(\mathcal{M})$ is
metrizable.
\end{corollary}

Sikora noticed in \cite{Sikora} that Corollary \ref{Metr} holds in the
special case when $\mathcal{M}$ is a countable semigroup.\bigskip 

\begin{proof}%
Since $LO(\mathcal{M})$ embeds into the Cantor set, $LO(\mathcal{M})$ is
metrizable.%
\end{proof}%
\bigskip

We conclude this section with the following lemma, which we use in Section
5. First we describe the notions of direct product and direct sum of magmas.
Let $\Gamma $ be a well-ordered set of indices, and let $\{\mathcal{M}%
_{\alpha }\}_{\alpha \in \Gamma }$ be a family of non-empty magmas. We
denote by $\prod\nolimits_{\alpha \in \Gamma }\mathcal{M}_{\alpha }$ the
direct product of magmas. For each $\alpha \in \Gamma $ we choose $b_{\alpha
}\in \mathcal{M}_{\alpha }$. Then we consider the direct sum of magmas $%
\bigoplus\nolimits_{\alpha \in \Gamma }\mathcal{M}_{\alpha }$ to be the
submagma of the magma $\prod\nolimits_{\alpha \in \Gamma }\mathcal{M}%
_{\alpha }$ generated by the sequences $\{x_{\alpha }\}_{\alpha \in \Gamma }$
where all but finitely many $x_{\alpha }$'s are in the submagmas of $\mathcal{%
M}_{\alpha }$ generated by $b_{\alpha }$. In the case of semigroups with
identity we take $b_{\alpha }$ to be the corresponding identity element. In
examples discussed in the next section (quandles) we always have $xx=x$.
Therefore, the submagma generated by $b_{\alpha }$ is a one-element
submagma. For each $\alpha \in \Gamma $, let $\mathcal{R}_{\alpha }\in LO(%
\mathcal{M}_{\alpha })$. Define the lexicographic order $\mathcal{R}$ on $%
\prod\nolimits_{\alpha \in \Gamma }\mathcal{M}_{\alpha }$ as follows. For $%
\{x_{\alpha }\}_{\alpha \in \Gamma },\{y_{\alpha }\}_{\alpha \in \Gamma }\in
\prod\nolimits_{\alpha \in \Gamma }\mathcal{M}_{\alpha }$ we say $%
(\{x_{\alpha }\}_{\alpha \in \Gamma },\{y_{\alpha }\}_{\alpha \in \Gamma
})\in \mathcal{R}$ if and only if for the smallest $\alpha \in \Gamma $ for
which $x_{\alpha }\neq y_{\alpha }$ we have $(x_{\alpha },y_{\alpha })\in 
\mathcal{R}_{\alpha }$. We denote the space of lexicographic left strict
orders on $\prod\nolimits_{\alpha \in \Gamma }\mathcal{M}_{\alpha }$ by $%
LO^{lex}(\prod\nolimits_{\alpha \in \Gamma }\mathcal{M}_{\alpha })$.
Analogously, we define the subspace $LO^{lex}(\bigoplus\nolimits_{\alpha \in
\Gamma }\mathcal{M}_{\alpha })$ of lexicographic left strict orders on $%
\bigoplus\nolimits_{\alpha \in \Gamma }\mathcal{M}_{\alpha }$, and we notice
that both spaces are homeomorphic. The homeomorphism is given by restricting
lexicographical orders on the direct product to the direct sum.

\begin{lemma}
\label{magma_products}$\mathbf{(i)}$ The product $\prod\nolimits_{\alpha \in
\Gamma }LO(\mathcal{M}_{\alpha })$ with Tychonoff topology is a closed
subspace of $LO(\bigoplus\nolimits_{\alpha \in \Gamma }\mathcal{M}_{\alpha
}) $.

$\mathbf{(ii)}$ The product $\prod\nolimits_{\alpha \in \Gamma }LO(\mathcal{M%
}_{\alpha })$ with Tychonoff topology is a closed subspace of $%
LO(\prod\nolimits_{\alpha \in \Gamma }\mathcal{M}_{\alpha })$.
\end{lemma}

\begin{proof}%
(i) We show that the subspace $LO^{lex}(\bigoplus\nolimits_{\alpha \in
\Gamma }\mathcal{M}_{\alpha })$ of the space $LO(\bigoplus\nolimits_{\alpha
\in \Gamma }\mathcal{M}_{\alpha })$ is homeomorphic to $\prod\nolimits_{%
\alpha \in \Gamma }LO(\mathcal{M}_{\alpha })$. The homeomorphism%
\begin{equation*}
\Psi :\prod\nolimits_{\alpha \in \Gamma }LO(\mathcal{M}_{\alpha
})\rightarrow LO^{lex}(\bigoplus\nolimits_{\alpha \in \Gamma }\mathcal{M}%
_{\alpha })
\end{equation*}%
is given by 
\begin{equation*}
\Psi (\{\mathcal{R}_{\alpha }\}_{\alpha \in \Gamma })\text{ is the
lexicographic order corresponding to }\{\mathcal{R}_{\alpha }\}_{\alpha \in
\Gamma. }
\end{equation*}%
Bijection of $\Psi$ follows directly from the definition, 
and the continuity of $\Psi $
follows from the fact that the image of the subbasis element $%
S_{(a,b)}^{\alpha }$ determined by $S_{(a,b)}\subset LO(\mathcal{M}_{\alpha
})$ is open. To see this, notice that if $\mathcal{R\in }\Psi
(S_{(a,b)}^{\beta })\subset LO(\prod\nolimits_{\alpha \in \Gamma }\mathcal{M}%
_{\alpha })$, then the subbasis element $U_{\mathcal{R}}=S_{(\{x_{\alpha
}\}_{\alpha \in \Gamma },\{y_{\alpha }\}_{\alpha \in \Gamma })}$, where $%
x_{\alpha }=y_{\alpha }=b_{\alpha }$ ($\alpha \neq \beta $) and $x_{\beta
}=a $, $y_{\beta }=b$, is a clopen neighborhood of $\mathcal{R}$ in $%
LO(\prod\nolimits_{\alpha \in \Gamma }\mathcal{M}_{\alpha })$ whose
restriction to $LO^{lex}(\bigoplus\nolimits_{\alpha \in \Gamma }\mathcal{M}%
_{\alpha })$ is in $\Psi (S_{(a,b)}^{\beta })$.

(ii) This can be established by a similar argument as (i). 
\end{proof}%

\begin{corollary}
If for all $\alpha \in \Gamma $, we have $LO(\mathcal{M}_{\alpha })\neq
\emptyset $, then the Cantor cube $\{0,1\}^{\left\vert \Gamma \right\vert }$
is a subset of $LO(\bigoplus\nolimits_{\alpha \in \Gamma }\mathcal{M}%
_{\alpha })$.
\end{corollary}

\begin{proof}%
Since $\left\vert LO(\mathcal{M}_{\alpha })\right\vert \geq 2$, we have $%
\{0,1\}^{\left\vert \Gamma \right\vert }\subset \prod\nolimits_{\alpha \in
\Gamma }LO(\mathcal{M}_{\alpha })$, and then we use Lemma \ref%
{magma_products}(i).%
\end{proof}%

\section{Quandles}

Important examples of magmas come from knot theory where they are used to
produce invariants of links. They are known as quandles and were introduced
and first studied by Joyce \cite{Joyce} and Matveev \cite{Matveev}. Recall
that a set $\mathcal{Q}$ with a binary operation $\ast :\mathcal{Q}\times 
\mathcal{Q}\rightarrow \mathcal{Q}$ is called a \emph{quandle} if the
following conditions are satisfied:

\begin{description}
\item[(i)] $(\forall a)\mathcal{[}a\ast a=a]$ (idempotence property);

\item[(ii)] For every $b\in \mathcal{Q}$, the mapping $\ast _{b}:\mathcal{%
Q\rightarrow Q}$ defined by $\ast _{b}(a)=a\ast b$ is bijective;

\item[(iii)] $(\forall a,b,c\mathcal{)[}(a\ast b)\ast c=(a\ast c)\ast (b\ast
c)]$ (right self-distributivity).
\end{description}

We denote the inverse of $\ast _{b}$ by $\overline{\ast }_{b}$ and write $%
\overline{\ast }_{b}(a)=a\overline{\ast }b$.\bigskip

\noindent \textbf{Remark 1 }The results of this section could also be
applied to more general structures of the so-called \emph{racks} 
(when axiom $(
\mathbf{i})$ is omitted in the above definition of a quandle). Racks were
defined by Fenn and Rourke in \cite{Fenn-Rourke}.\bigskip 

\noindent \textbf{Example} \textbf{1 }A quandle $\mathcal{Q}$ is called 
\emph{trivial} if the quandle operation $\ast $ is defined by $(\forall
a,b)[a\ast b=a]$. Then every total strict order on $\mathcal{Q}$ is right
invariant under the quandle operation. Moreover, $RO(\mathcal{Q})=\{0,1\}^{|%
\mathcal{Q}|}$ for $|\mathcal{Q}|=\aleph _{0}$, because in this case $RO(%
\mathcal{Q})$ is a zero-dimensional, compact, separable topological space
without isolated points.\bigskip

Let $G$ be a group. We define the \emph{conjugate} quandle $\mathrm{Conj}(G)$
as one with domain $G$ and the quandle operation $\ast $ given by $a\ast
b=b^{-1}ab$.

\begin{proposition}
Let $G$ be a bi-orderable group. Then $\mathrm{Conj}(G)$ is right orderable
and every bi-order on $G$ induces a right order on $\mathrm{Conj}(G)$.
\end{proposition}

\begin{proof}%
Let $\mathcal{P}$ be a bi-order on $G$. Then, by the definition of $\mathcal{%
P}$, we have%
\begin{equation*}
(\forall a,b,c)[(a,b)\in \mathcal{P}\Rightarrow (c^{-1}ac,c^{-1}bc)\in 
\mathcal{P})]
\end{equation*}%
Using $\mathcal{P}$, we define $\mathcal{R}$ on $\mathrm{Conj}(G)$ by
\begin{equation*}
(\forall a,b)[(a,b)\in \mathcal{R\Leftrightarrow (}e,a^{-1}b)\in \mathcal{P}%
],
\end{equation*}%
where $e$ is the identity of $G$. The order $\mathcal{R}$ is right invariant
because for $(a,b)\in \mathcal{R}$ and $c\in \mathrm{Conj}(G)$, we have 
\begin{equation*}
(e,(a\ast c)^{-1}(b\ast
c))=(e,(c^{-1}a^{-1}c)(c^{-1}bc))=(e,c^{-1}(a^{-1}b)c)\in \mathcal{P}\text{.}
\end{equation*}%
Since $(e,a^{-1}b)\in \mathcal{P}$, we have $(a\ast c,b\ast c)\in \mathcal{R}
$.%
\end{proof}%
\bigskip 

\noindent \textbf{Remark 2} Notice that not all right orders on $\mathrm{Conj%
}(G)$ are induced by bi-orders on $G$. It is possible to have $%
BiO(G)=\emptyset $ while $RO(\mathrm{Conj}(G))\neq \emptyset $. For example,
let $G$ be an abelian group with torsion. Then $BiO(G)=\emptyset ,$ but $%
\mathrm{Conj}(G)$ is a trivial quandle, so it admits many right orders (see
Example 1).\bigskip

\noindent \textbf{Example 2 }($n$-quandle \cite{Joyce}) Let $n$ be a
positive integer. Consider a quandle $\mathcal{Q}$ satisfying the following
identity for all $a,b$:\ $b\ast a^{\ast n}=b$, where $b\ast a^{\ast
n}=(...((b\ast \underset{n}{\underbrace{a)\ast a)\ast ...\ast a)\ast a}}$.
Then $RO(\mathcal{Q})=\emptyset $ unless $n=1$, in which case the quandle is
trivial. The case when $n=2$ plays an important role in knot theory, and the
quandle is called \emph{involutive} or \emph{Kei}. In particular, for every
group we define Kei by setting $b\ast a=ab^{-1}a$.\bigskip 

Example 2 motivates the following general result.

\begin{proposition}
\label{Not_ord}$\mathbf{(i)}$ Let $\mathcal{Q}$ be a right-orderable
quandle. If $b\ast a^{\ast n}=b$ for some elements $a,b\in \mathcal{Q}$ and
a positive integer $n$, then $b\ast a=b$.

$\mathbf{(ii)}$ Let $a$ and $b$ be elements of a group $G$ such that $a$ is
not in the centralizer $C_{G}(b)$ of $b$, but $a^{n}\in C_{G}(b)$ for some $%
n>1$. Then \textrm{Conj}$(G)$ is not right orderable.
\end{proposition}

\begin{proof}%
(i) Let $\mathcal{R}$ be a strict right order on $\mathcal{Q}$. If $b\ast
a\neq b$, then either $(b,b\ast a)\in \mathcal{R}$ or $(b\ast a,b)\in 
\mathcal{R}$. Consider the case when $(b,b\ast a)\in \mathcal{R}$. Then, by
right invariance, we have%
\begin{equation*}
(b\ast a,b\ast a^{\ast 2})\in \mathcal{R\wedge }(b\ast a^{\ast 2},b\ast
a^{\ast 3})\in \mathcal{R\wedge }...\wedge (b\ast a^{\ast (n-1)},b\ast
a^{\ast n})\in \mathcal{R}\text{.}
\end{equation*}%
Hence, by transitivity, $(b,b\ast a^{\ast n})=(b,b)\in \mathcal{R}$, which
is a contradiction. Similarly, if $(b\ast a,b)\in \mathcal{R}$, then $(b\ast
a^{\ast n},b)=(b,b)\in \mathcal{R}$, which is a contradiction. Thus $b\ast
a=b$.

(ii) Since $a\in G\backslash C_{G}(b)$, there is $b\in G$ such that $ba\neq
ab$ (i.e., $b\ast a\neq b$). Since $a^{n}\in C_{G}(a)$, we have $b\ast
a^{\ast n}=a^{-n}ba^{n}=b$, and, by (i), $\mathrm{Conj}(G)$ is not right
orderable.%
\end{proof}%
\bigskip 

B. H. Neumann noted in \cite{Neumann} (Lemma 1.1) that if $G$ is a
bi-orderable group, $a,b\in G$, and $b$ commutes with $a^{n}$ ($n\neq 0$),
then $b$ commutes with $a$. This also follows from Proposition \ref{Not_ord}%
(ii). Namely, if $G$ is bi-orderable, then $\mathrm{Conj}(G)$ is right
orderable. The condition that $a^{n}$ commutes with $b$ is equivalent to $%
b\ast a^{\ast n}=b$, hence $n=1$ and we obtain $b=b\ast a=a^{-1}ba$, so $b$
commutes with $a$. Neumann also observed (\cite{Neumann}, Lemma 1.3) that if 
$G$ is a bi-orderable group and $a$ commutes with the commutator $%
[a^{n},b]=a^{-n}b^{-1}a^{n}b$ then $a$ commutes with $[a,b]$. This is
also a consequence of Proposition \ref{Not_ord}$(ii)$ applied to $(x\ast
y)\ast x^{\ast n}$ for $b=x\ast y$ and $a=x$. Finally, Neumann asked whether
in a bi-orderable group $G,$ if $a$ commutes with $[[a^{n},b],a]$ then $a$
commutes with $[[a,b],a]$. This question (and its generalized version) was
answered negatively by Mura and Rhemtulla (see \cite{MR}, Lemma 2.5.3). In
the language of quandles this shows, in particular, that there is a right
orderable quandle $\mathrm{Conj}(G)$ in which the identity $(((a\overline{%
\ast }b)\overline{\ast }a^{\ast n})\ast b)\ast a=((a\overline{\ast }b)%
\overline{\ast }a^{\ast n})\ast b$ does not imply $(((a\overline{\ast }b)%
\overline{\ast }a)\ast b)\ast a=((a\overline{\ast }b)\overline{\ast }a)\ast b
$.\bigskip 

\noindent \textbf{Remark 3} As we observed in Remark 2, every bi-order on a
group $G$ induces a right invariant order on \textrm{Conj}$(G)$. However,
the existence of left orders on $G$ is not sufficient for the existence of
right orders on $\mathrm{Conj}(G)$. For example, consider the fundamental
group of the Klein bottle%
\begin{equation*}
Kb=\left\langle a,b|\text{ }a^{-1}ba=b^{-1}\right\rangle \text{,}
\end{equation*}%
which is left orderable, but not bi-orderable. By Proposition \ref{Not_ord}%
(ii), $RO(\mathrm{Conj}(Kb))=\emptyset $ because $a^{-2}ba^{2}=b$ or, using
quandle operation, $b\ast a^{\ast 2}=b$.\bigskip

\noindent \textbf{Example 3} Another application of Proposition \ref{Not_ord}%
(ii) is the quandle $\mathrm{Conj}(G_{(n,m)})$ for $n,m>1,$ where 
\begin{equation*}
G_{(n,m)}=\left\langle x,y|\text{ }x^{n}=y^{m}\right\rangle \text{.}
\end{equation*}%
Since $x\notin Z(G_{(n,m)})$ but $x^{n}\in Z(G_{(n,m)})$, it follows that $%
RO(\mathrm{Conj}(G_{(n,m)}))=\emptyset $. In the case when $n=m=2$ we obtain
the fundamental group of the Klein bottle. When $n$ and $m$ are relatively
prime, the group $G_{(n,m)}$ is the fundamental group of $(n,m)$-torus knot.
In the case when $n=2$ and $m=3$ we have Artin's braid group $B_{3}$. Since $%
B_{3}$ is a subgroup of the braid group $B_{n}$ for $n\geq 3$, \textrm{Conj}(%
$B_{n})$ is not right orderable. Recall that $B_{n}$ is right orderable but
not bi-orderable \cite{Dynnikov}.

\section{Compactness and Conrad's theorem}

In Theorem \ref{Embedding} we proved that the space of left orders on a
magma is compact. In the case of a group $G$ we found an interesting
connection between compactness of the space $LO(G)$ and the classical
theorem due to Conrad and \L o\'{s} (\cite{Conrad}, \cite{Los}; see Theorem %
\ref{Conrad} below). Let $P\subset G$ be a sub-semigroup that is \emph{pure}
(i.e., if $g\in P$ then $g^{-1}\notin P$). Such a subset $P$ is the \emph{%
positive} \emph{cone} of a strict partial left order $<$ on $G$ determined
by 
\begin{equation*}
x<y\Leftrightarrow x^{-1}y\in P.
\end{equation*}%
Notice that if $x<y$, then for all $z$, we have $zx<zy$. Conversely, every
partial left order on $G$ determines a positive cone. A positive cone $P^{+}$
is \emph{total} if $P^{+}$ contains $g$ or $g^{-1}$ for each non-identity
element $g\in G$. A pure and positive cone uniquely determines a strict
linear left order on $G$. For a subset $A\subset G$, we denote by $sgr(A)$
the sub-semigroup of $G$ generated by $A$.

\begin{theorem}
\label{Conrad}$($Conrad $)$ A partial left order $P$ on $G$\ can be extended
to a total left order if and only if for every finite set $%
\{x_{1},x_{2},...,x_{n}\}\subset G\backslash \{e\}$, there is a
corresponding sequence $\epsilon _{1},\epsilon _{2},...,\epsilon _{n}$ $%
(\epsilon _{i}\in \{+1,-1\},$ $i=1,2,...,n)$ such that 
\begin{equation*}
e\notin sgr(P\cup \{x_{1}^{\epsilon _{1}},x_{2}^{\epsilon
_{2}},...,x_{n}^{\epsilon _{n}}\}).
\end{equation*}
\end{theorem}

Not every partial left order on $G$ extends to a linear left order on $G$,
even if $G$ admits a left invariant order. For example, $P=sgr(%
\{b^{2},a,ab^{-2}\})\subset $ $Kb$ is a cone of a strict partial left order
on $G$, which cannot be extended to a linear left order on $G$.

We observe that the algebraic condition given in Theorem \ref{Conrad}
resembles the condition for compactness in terms of closed sets: A Hausdorff
topological space $X$ is compact if every family of its closed subsets with
the finite intersection property (i.e., every finite subfamily has 
a non-empty
intersection) has a non-empty intersection. In fact, we can use Conrad's
theorem to show that $LO(G)$ is compact. We sketch a proof below.

By Alexander's Subbase Theorem \cite{Alex} we can reduce the question about
compactness of $LO(G)$ to the families of subsets of $LO(G)$ consisting of
elements of a subbasis: 
\begin{equation*}
\mathcal{A=\{}S_{g_{\alpha }}\mathcal{\}}_{\alpha \in \Gamma },\text{ where }%
S_{g_{\alpha }}=\{P^{+}\in LO(G)|\text{ }g_{\alpha }\in P^{+}\}\in \mathcal{S%
}\text{ for }\alpha \in \Gamma ,
\end{equation*}%
(which are clopen subsets). We notice that the finite intersection property
for the family $\mathcal{A}$ allows us to show that the algebraic condition
given in Conrad's theorem holds for $P=sgr(\{g_{\alpha }|$ $\alpha \in
\Gamma \})$. Namely, by the finite intersection property, for any finite set
of elements $g_{\alpha _{1}},g_{\alpha _{2}},...,g_{\alpha _{k}}$ ($\alpha
_{i}\in \Gamma $), there is a strict linear left order $Q^{+}$ of $G$ which,
contains $sgr(\{g_{\alpha _{1}},g_{\alpha _{2}},...,g_{\alpha _{k}}\})$. If $%
\{x_{1},x_{2},...,x_{n}\}\subset G\backslash \{e\}$, then there is a
corresponding sequence $\epsilon _{1},\epsilon _{2},...,\epsilon _{n}$ $%
(\epsilon _{i}\in \{+1,-1\},$ $i=1,2,...,n)$ such that $x_{i}^{\epsilon
_{i}}\in Q^{+}$. Therefore, the condition $e\notin sgr(P\cup
\{x_{1}^{\epsilon _{1}},x_{2}^{\epsilon _{2}},...,x_{n}^{\epsilon _{n}}\})$
holds for $P$. By Theorem \ref{Conrad},
\begin{equation*}
\bigcap\limits_{\alpha \in \Gamma }S_{g_{\alpha }}\neq \emptyset .
\end{equation*}

\section{Examples and open problems}

In some cases, we are able to provide a complete characterization of the
spaces of strict left orders. In particular, this is true in the case when $%
LO(\mathcal{M})$ (or $RO(\mathcal{M})$) has a countable basis and no
isolated points. In this case the set of left (right) orders is homeomorphic
to the Cantor set $\{0,1\}^{\aleph _{0}}$. We have illustrated this in
Example 1 when $\mathcal{M}$ is an infinite countable trivial quandle. Let $%
\mathbb{Z}
^{\oplus \mathfrak{m}}=\oplus _{\mathfrak{m}}%
\mathbb{Z}
$ denote the free abelian group with basis of cardinality $\mathfrak{m}$. If 
$1<\mathfrak{m\leq \aleph }_{0}$ then the space of bi-orders $BiO(%
\mathbb{Z}
^{\oplus \mathfrak{m}})$ is homeomorphic to the Cantor set (\cite{Sikora}, 
\cite{Dabkowska}). We also note that $BiO(%
\mathbb{Q}
^{\oplus \mathfrak{m}})$ is homeomorphic to $BiO(%
\mathbb{Z}
^{\oplus \mathfrak{m}})$ (see \cite{KopMed}). This observation has the
following proposition as a consequence.

\begin{proposition}
$\mathbf{(i)}$ The space $BiO(%
\mathbb{Q}
^{\oplus \mathfrak{m}})$ is metrizable if and only if $\mathfrak{m\leq
\aleph }_{0}$.

$\mathbf{(ii)}$ The space $BiO(%
\mathbb{R}
)$ is not metrizable.
\end{proposition}

\begin{proof}
(i) If $\mathfrak{m>\aleph }_{0}$, we use Lemma \ref{magma_products}(i) for $%
\mathcal{M}_{\alpha }=%
\mathbb{Q}
$ and $BiO(%
\mathbb{Q}
)=\{0,1\}$, that is, $\{0,1\}^{\mathfrak{m}}=BiO^{lex}(%
\mathbb{Q}
^{\oplus \mathfrak{m}})$ $\subset BiO(%
\mathbb{Q}
^{\oplus \mathfrak{m}})$. Since the Cantor cube $\{0,1\}^{\mathfrak{m}}$ is
not metrizable, neither is $BiO(%
\mathbb{Q}
^{\oplus \mathfrak{m}})$.

(ii) This statement follows from the fact that $%
\mathbb{R}
$ is isomorphic to $%
\mathbb{Q}
^{\oplus \mathfrak{c}}$, where $\mathfrak{c=}\left\vert 
\mathbb{R}
\right\vert $.%
\end{proof}%
\bigskip 

We showed that for $\mathfrak{m\geq \aleph }_{0}$, we have the following 
embeddings 
\begin{equation*}
\{0,1\}^{\mathfrak{m}}\hookrightarrow BiO(%
\mathbb{Z}
^{\oplus \mathfrak{m}})\hookrightarrow \{0,1\}^{\mathfrak{m}}\text{.}
\end{equation*}%
Furthermore, $BiO(%
\mathbb{Z}
^{\oplus \mathfrak{m}})$ has no isolated points. For $\mathfrak{m=\aleph }%
_{0}$ this implies that $BiO(%
\mathbb{Z}
^{\oplus \mathfrak{m}})$ is homeomorphic to the Cantor set. In general, we
propose the following problem.

\begin{problem}
Is $LO(%
\mathbb{Z}
^{\oplus \mathfrak{m}})$ homeomorphic to $\{0,1\}^{\mathfrak{m}}$ for all
cardinals $\mathfrak{m}>\aleph _{0}$?
\end{problem}

Recall that a Hausdorff topological space $X$ is called \emph{supercompact}
if there exists a subbasis $\mathcal{S}$ such that for each covering $%
\mathcal{A}$ consisting of elements of $\mathcal{S}$, there exists a
subcovering consisting of two elements. The Cantor cube is supercompact, and 
$LO(\mathcal{M})$ is supercompact provided it has a countable weight (it
follows from \cite{Szymanski} and \cite{Mill} that every compact metric
space is supercompact).

\begin{problem}
For which magma $\mathcal{M}$ is the space $LO(\mathcal{M})$ supercompact$?$%
\bigskip
\end{problem}

\noindent \textbf{Acknowledgement} Harizanov was partially supported by the
NSF grant 0502499, and by the UFF\ award of the George Washington University.

\end{document}